# COMPARING HEEGAARD AND JSJ STRUCTURES OF ORIENTABLE 3-MANIFOLDS

MARTIN SCHARLEMANN AND JENNIFER SCHULTENS

ABSTRACT. The Heegaard genus $g$ of an irreducible closed orientable 3-manifold puts a limit on the number and complexity of the pieces that arise in the Jaco-Shalen-Johannson decomposition of the manifold by its canonical tori. For example, if $p$ of the complementary components are not Seifert fibered, then $p \leq g - 1$. This generalizes work of Kobayashi [Ko]. The Heegaard genus $g$ also puts explicit bounds on the complexity of the Seifert pieces. For example, if the union of the Seifert pieces has base space $P$ and $f$ exceptional fibers, then
$$f - \chi(P) \leq 3g - 3 - p.$$

## 1. INTRODUCTION

Nearly a century ago, Heegaard noticed that quite a few 3-manifolds could be written as the union of two handlebodies ([He], see also [Prz] for a translation of the relevant parts) . Later it was discovered that this first global structure theorem applied in fact to all 3-manifolds. Now called a *Heegaard splitting* of the 3-manifold, this structure has proven to be a deceptively simple picture because, although the existence of the structure is easy to prove, it is not unique. A single manifold may have several Heegaard splittings and the relationship between the various splittings has been difficult to understand.

A modern and more useful structure theorem for 3-manifolds, due to Jaco-Shalen and Johannson, does not have the problem of non-uniqueness. In its simplest form the theorem states that, for any irreducible orientable closed 3-manifold $M$ there is a collection $\Theta$ of incompressible tori (called the *canonical tori* of $M$) so that each complementary component of $M - \Theta$ is either a Seifert manifold (possibly without exceptional fibers) or is both acylindrical (any properly imbedded incompressible annulus is boundary parallel) and atoroidal (any properly imbedded torus is boundary parallel). Moreover, $\Theta$ is unique up to ambient isotopy in $M$.

The connection between Heegaard structure and this "JSJ structure" has been poorly understood. The only significant information comes from two theorems of Kobayashi [Ko, Theorems 1 and 2]. The first states that if a closed orientable 3-manifold has a genus $g$ Heegaard splitting then $\Theta$ has at most $3g - 3$ complementary components. Moreover, if it has exactly $3g-3$, then every complementary component is atoroidal (though not necessarily acylindrical); the second theorem gives more detail about the structure of these $3g - 3$ components, particularly the non-Seifert pieces.

Research supported in part by NSF grants and MSRI.

In part, we here expand on Kobayashi's theme. For example, we show that if $M$ has a genus $g$ Heegaard splitting, then at most $g-1$ complementary components are not Seifert manifolds. Moreover, among the Seifert pieces, if $n'$ fiber over the twice punctured projective plane or, with one exceptional fiber, over the once punctured projective plane, and $n$ other components are also toroidal, then the number of complementary components is no more than $3g - 3 - n - n'/2$. (See Corollary 6.6.) We know no counterexamples to the stronger statement that there are no more than $3g - 3 - n - n'$ toroidal complementary components.

In addition, we find limits on the complexity of the Seifert pieces. In order of increasing generality, these limits occur as Theorems 4.7, 5.3, and 6.4. Mostly to obtain the limits on the number of toroidal components mentioned above, the results are a bit stronger than the following more easily stated corollary: For $p$ the number of non-Seifert components, $P$ the base surface of the Seifert parts, and $f$ the number of exceptional fibers, we have $f - \chi(P) \leq 3g - 3 - p$. (See Corollary 6.9.)

Finally, Kobayashi's structure theorem ([Ko, Theorem 2]) for "full" Haken manifolds is shown to have the following beautiful explanation (see Corollary 4.8): When an irreducible Heegaard splitting of a full Haken 3-manifold is put in thin position (see [ST]), then it is strung out like an array of jewelry: The setting consists of Seifert pieces connected together by amalgamating tori. Embedded in this setting are $g-1$ "jewels", each homeomorphic to the complement of a 2-bridge link in $S^3$.

Here is an outline: In Section 2 we briefly recount the theory of generalized Heegaard splittings and untelescopings, mostly from [ST]. The first of two core sections is Section 3, wherein we explain the delicate process of positioning the canonical tori optimally with respect to the surfaces that arise from a strongly irreducible generalized Heegaard splitting. In the end, the Seifert pieces can have one of two positions, aligned or non-aligned. Section 4 explains the connection between the complexity of Seifert pieces and how they intersect the compression-bodies of the generalized splitting. The results here apply only to the aligned pieces, but they are the critical ones. The remainder of the paper is designed to incorporate the non-aligned pieces so that ultimately we are able to make statements that do not require any knowledge about which Seifert pieces are aligned and which are not. So, as the argument progresses, the statements get a bit weaker but the generality with which they can be applied improves. By the end of Section 6 the distinction between the aligned and non-aligned positions of the Seifert pieces no longer needs to be considered; an inequality shows (roughly) that the vertical index, with which the complexity of the aligned pieces is measured, bounds the horizontal index, with which the complexity of the non-aligned pieces is measured.

## 2. Heegaard splittings and their untelescopings

**Definition 2.1.** *A compression body $H$ is a connected 3-manifold obtained from a closed surface $\partial_- H$ by attaching 1-handles to $\partial_- H \times \{1\} \subset \partial_- H \times I$. Dually, a compression body is obtained from a connected surface $\partial_+ H$ by attaching 2-handles to $\partial_+ H \times \{1\} \subset \partial_+ H \times I$ and 3-handles to any 2-spheres thereby created. The cores of the 2-handles are called* meridian disks

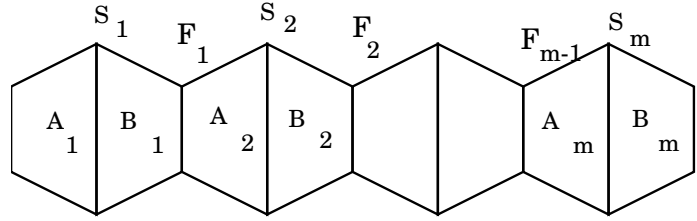

FIGURE 1.

*For $H$ a compression body, define the index $J(H) = \chi(\partial_- H) - \chi(\partial_+ H)$.*

A *Heegaard splitting* $M = A \cup_S B$ of a compact orientable 3-manifold consists of an orientable surface $S$ in $M$, together with two compression bodies $A$ and $B$ so that $S = \partial_+ A = \partial_+ B$ and $M = A \cup_S B$. $S$ itself is called the splitting surface. The genus of the splitting is defined to be the genus of $S$.

A *stabilization* of $A \cup_S B$ is the Heegaard splitting obtained by adding to $A$ a regular neighborhood of a proper arc in $B$ which is parallel in $B$ to an arc in $S$. A stabilization has genus one larger and, up to isotopy, is independent of the choice of arc in $B$. If the construction is done symmetrically to an arc in $A$ instead, the two splittings are isotopic.

Recall the following (see e. g. [Sc]): If there are meridian disks $D_A$ and $D_B$ in $A$ and $B$ respectively so that $\partial D_A$ and $\partial D_B$ intersect in a single point in $S$, then $A \cup_S B$ can be obtained by stabilizing a lower genus Heegaard splitting. We then say that $A \cup_S B$ is *stabilized*. If there are meridian disks $D_A$ and $D_B$ in $A$ and $B$ respectively so that $\partial D_A$ and $\partial D_B$ are disjoint in $S$, then (see [CG]) $A \cup_S B$ is *weakly reducible*. If there are meridian disks so that $\partial D_A = \partial D_B$, then $A \cup_S B$ is *reducible*. It is easy to see that reducible splittings are weakly reducible and that (except for the genus one splitting of $S^3$) any stabilized splitting is reducible. It is a theorem of Haken [Ha] that any Heegaard splitting of a reducible 3-manifold is reducible and it follows from a theorem of Waldhausen [W] that a reducible splitting of an irreducible manifold is stabilized.

**Definition 2.2.** *Suppose $M$ is an irreducible closed orientable 3-manifold. A generalized Heegaard splitting of $M$ is a structure*

$$M = (A_1 \cup_{S_1} B_1) \cup_{F_1} (A_2 \cup_{S_2} B_2) \cup_{F_2} ... \cup_{F_{m-1}} (A_m \cup_{S_m} B_m).$$

*Here each $A_i$ and $B_i$ is a compression body, $\partial_+ A_i = S_i = \partial_+ B_i$ (so $A_i \cup_{S_i} B_i$ is a Heegaard splitting), $\partial_- B_i = F_i = \partial_- A_{i+1}$, and each $F_i$ is incompresible in $M$. We say that a generalized Heegaard splitting is* strongly irreducible *if each Heegaard splitting $A_i \cup_{S_i} B_i$ is strongly irreducible. For a generalized splitting we will often denote $S = \cup_{i=1}^m S_i$ and $F = \cup_{i=1}^{m-1} F_i$. (See Figure 1.)*

The central theorem of [ST], via the calculation of [Sc3, Lemma 2], directly implies this:

**Theorem 2.3.** *Suppose $M$ is an irreducible closed orientable 3-manifold and $M$ has a genus $g$ Heegaard splitting. Then $M$ has a strongly irreducible generalized Heegaard*

*splitting*

$$(A_1 \cup_{S_1} B_1) \cup_{F_1} (A_2 \cup_{S_2} B_2) \cup_{F_2} ... \cup_{F_{m-1}} (A_m \cup_{S_m} B_m)$$

*so that*

$$\sum_{i=1}^{m} J(A_i) = \sum_{i=1}^{m} J(B_i) = \chi(F) - \chi(S) \leq 2g - 2.$$

This structure is created by *untelescoping* a minimal genus Heegaard splitting of $M$. Notice that the theorem is a tautology if $M$ has a genus $\leq g$ Heegaard splitting that is strongly irreducible.

## 3. Heegaard splitting surfaces vs. canonical tori

We would like to understand how the surfaces $S$ and $F$ of a generalized strongly irreducible Heegaard splitting of $M$ intersect the canonical tori $\Theta$ of $M$. More generally, we would like to simplify as much as possible the intersections of $S$ and $F$ with those complementary components of $\Theta$ in $M$ that are Seifert.

We begin with a fairly easy argument:

**Definition 3.1.** *A properly imbedded collection of annuli in a 3-manifold $M$ is essential if it is incompressible and no component is $\partial$-parallel.*

**Lemma 3.2.** *Suppose $M$ is an irreducible closed orientable 3-manifold with a strongly irreducible generalized Heegaard splitting*

$$(A_1 \cup_{S_1} B_1) \cup_{F_1} (A_2 \cup_{S_2} B_2) \cup_{F_2} ... \cup_{F_{m-1}} (A_m \cup_{S_m} B_m)$$

*and suppose $\Theta$ is the set of canonical tori of $M$. Then $\Theta$ may be isotoped so that $F \cap \Theta$ consists of curves essential in both $F$ and $\Theta$ and so that $\Theta$ intersects each compression body $A_i$ and $B_i$ only in essential annuli and incompressible tori.*

Note that an incompressible torus in a compression body $H$ must be parallel to a torus component of $\partial_- H$.

**Proof:** Here is a sketch. More detail can be found in e. g. [Sc4]. First note that, since both $F$ and $\Theta$ are incompressible, a simple innermost disk argument can be used to remove all components of $F \cap \Theta$ that are inessential in either surface (hence both surfaces). So we can assume that all components of $F \cap \Theta$ are essential in both $F$ and $\Theta$. Now the surface $S_i$ can be used to "sweep out" the region between $F_{i-1}$ and $F_i$, once certain 1-complexes incident to $F_{i-1}$ and $F_i$ are removed (the spines of $A_i$ and $B_i$ respectively). At the beginning of the sweep-out, each component of $S_i \cap \Theta$ is either essential or bounds a tiny disk in $A_i$, each disk corresponding to an intersection point of the spine of $A_i$ with $\Theta$. Similarly, at the end of the sweep-out, each component of $S_i \cap \Theta$ is either essential or bounds a tiny disk in $B_i$, each disk corresponding to an intersection point of the spine of $B_i$ with $\Theta$. There cannot simultaneously be disk components of intersection of $\Theta$ with $A_i$ and with $B_i$, since $A_i \cup_{S_i} B_i$ is strongly irreducible. So in some positioning, all components of intersection are essential.

Repeat this argument for each $S_i, 1 \leq i \leq m$. Then each component of $(F \cup S) \cap \Theta$ is essential in $\Theta$ and hence in $F \cup S$. At this point we know that $\Theta$ intersects each compression body $A_i$ and $B_i$ only in incompressible annuli and tori. Now remove any boundary parallel annuli by isotopies. □

**Definition 3.3.** *Suppose $V$ is a Seifert manifold with base space $P$. A surface $T \subset V$ is* vertical *if it is a union of generic fibers and is* horizontal *if it is transverse to each fiber.*

Note that a vertical surface then must have Euler characteristic zero, and so is a union of annuli and tori.

**Definition 3.4.** *Suppose $E$ is a possibly non-orientable surface and $\xi$ is an $I$-bundle over $E$ whose total space is orientable. Then let $\partial \xi$ denote the restriction of $\xi$ to $\partial E$ and let $\dot\xi$ denote the associated $\partial I$-bundle of $E$. The boundary of the total space of $\xi$ is the union of the two.*

**Theorem 3.5.** [Ja, Theorem VI.34] *Any properly imbedded incompressible and $\partial$-incompressible 2-sided surface in an orientable Seifert manifold $V$ can be properly isotoped so that either it is vertical or it is horizontal. If it is horizontal, then $V$ is the union, along $\dot\xi$, of two copies of an $I$-bundle $\xi$ over a surface $E$, and the incompressible surface consists of parallel copies of $\dot\xi$.*

Of course, if $\xi$ is an orientable $I$-bundle (so $E$ is orientable), then $V$ fibers over the circle with fiber $E$.

**Definition 3.6.** *When a Seifert manifold $V$ is expressed as the union, along $\dot\xi$, of two copies of an $I$-bundle $\xi$ over a surface $E$, we call $E$ the associated $I$-base of this construction, and say that $V$ is an $I$-bundle construct over $E$. $E$ is a branched cover of the Seifert base $P$ of $V$.*

For the purposes of this paper, we will always be able to assume that $\chi(E) \leq -1$, since if $E$ is the Möbius band or the annulus we could fiber $V$ differently so that $\dot\xi$ is vertical in $V$. Also, for expository purposes, little would be lost by always taking $\xi$ to be a product bundle, so that $(\xi; \dot\xi, \xi|\partial) = (E \times I; E \times \partial I, \partial E \times I)$.

**Corollary 3.7.** *Under the hypotheses of Lemma 3.2, $F$ may be further isotoped so that it intersects each Seifert component of $M - \Theta$ either in a horizontal surface or in vertical essential annuli and incompressible tori.*

**Proof:** Let $V$ be a Seifert component of $M - \Theta$ and isotope $\Theta$ as provided in Lemma 3.2. Then each component of $F \cap V$ is incompressible in $V$, since $F$ is incompressible in $M$ and all curves of $F \cap \Theta$ are essential in $F$. If any is $\partial$-compressible then, dually, some annulus $A$ in $\Theta - F$ $\partial$-compresses to $F$. The $\partial$-compression cannot turn $A$ into an essential disk, since $F$ is incompressible, so $A$ must be $\partial$-parallel in the submanifold $A_i \cup_{S_i} B_i$ in which it lies. It's easy to see that then either $A$ intersects either $A_i$ or $B_i$ in a $\partial$-parallel annulus or $A_i \cup_{S_i} B_i$ would be weakly reducible. Either case contradicts our hypotheses, so $A$ is also $\partial$-incompressible. The result then follows from Theorem 3.5. □

We would like now similarly to simplify the positioning of $S$ in each Seifert component of $M - \Theta$. We can come surprisingly close.

**Theorem 3.8.** *Assume the hypotheses of Lemma 3.7, and let $V$ be the union of Seifert components of $M - \Theta$. Then $S$ may be further isotoped so that for each $S_i$, every component but at most one of $S_i \cap V$ is either horizontal or vertical. The one exception is obtained either from one or two horizontal components by attaching a vertical tube between them, or from one or two vertical components by attaching a single horizontal tube, i. e. a tube whose core projects homeomorphically to an embedded arc in the base space $P$.*

This positioning of $\Theta$ with respect to the generalized Heegaard splitting is called a *preferred positioning*.

**Proof:** If $S_i \cap V$ is incompressible in $V$, the proof mimics that of Lemma 3.7. So suppose it is compressible and, with no loss of generality, suppose it compresses into $A_i \cap V$. After compressing maximally into $A_i \cap V$ it is standard to see (via strong irreducibility) that the resulting surface $\hat{S}_i$ is incompressible in $V$, so we can reconstruct $S_i$ by starting with an incompressible surface in $V$ and attaching tubes on one side (dual to the compressing disks we've just used.) Now consider $B_i$. By strong irreducibility, any meridian disk $D$ in $B_i$ must have a boundary that runs along each tube (i. e. intersects each compressing disk for $A_i$ in $V$). Either $D$ lies in $V$ or an outermost arc cuts off a $\partial$-compressing disk $D' \subset B_i$ which $\partial$-compresses an annulus of $\Theta - S_i$ to $S_i$ through $B_i$. Since the annulus cannot be $\partial$-parallel, the result is a compressing disk in $B_i$. If the disk were outside $V$, this would contradict strong irreducibility, so we deduce that it lies in $V$. The upshot is that we may assume that $S_i \cap V$ compresses in both $A_i \cap V$ and $B_i \cap V$ and, indeed, in the same component of $V$. So henceforth we may assume $V$ is connected.

Select a family of essential 2-sided arcs in the base space $P$ of $V$, chosen so that their complement consists entirely of disks, each containing at most one exceptional point and sufficiently plentiful that no disk lies on both sides of the same arc. Let $\mathcal{A} \subset V$ denote the family of vertical annuli that covers these arcs. The complementary components of $\mathcal{A}$ in $V$ are solid tori $T_1, ..., T_t$. Let $T = \cup_{j=1}^{t} T_j$.

**Claim:** $S_i$ can be isotoped so that it intersects each $\partial T_j$ only in essential curves. (See Figure 2.)

*Proof of claim:* This is an adaptation of the argument of Lemma 3.2. The details of the proof are a little more complex than the somewhat condensed version about to be given here. A more complicated argument in the same spirit, with full details, occurs in [RS, Section 2.2].

**Step 1:** Parameterize the sweep-out of Lemma 3.2 by the interval. Within it is a subinterval in which every curve in $S_i \cap \partial V$ is essential in both $S_i$ and $V$, but just before the subinterval a curve in $S_i \cap \partial V$ cuts off a meridian disk of $A_i$, say, lying in $\partial V$ and just after the subinterval a curve in $S_i \cap \partial V$ cuts off a meridian disk of $B_i$ lying in $\partial V$. Now consider how $S_i$ intersects $\mathcal{A}$ and $\cup_{j=1}^{t} \partial T_j$ during the sweep-out (now meaning sweep-out through this subinterval). There will be no critical points of intersection in $\partial V$ (since every component of $S_i \cap V$ remains essential in $\partial V$ throughout the sweep-out), so in fact we may as well take the curves $S_i \cap \partial V$ to be fixed throughout the sweep-out, and minimally intersecting the set of curves $\partial \mathcal{A}$. We can also remove by an isotopy all $\partial$-parallel annuli in $S_i \cap V$.

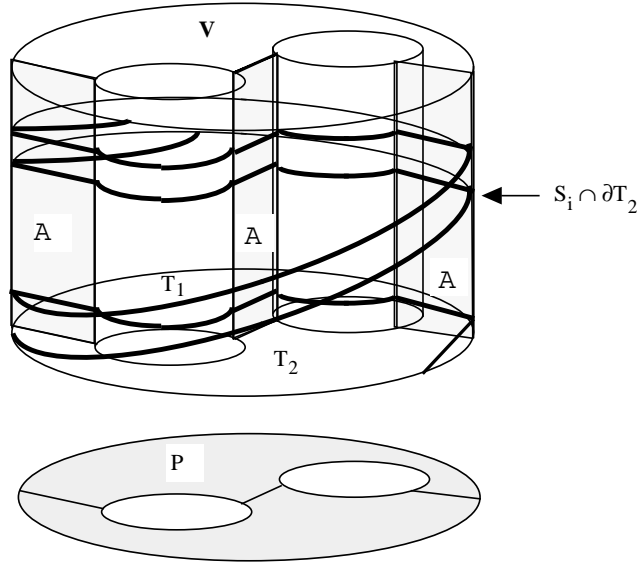

Figure 2.

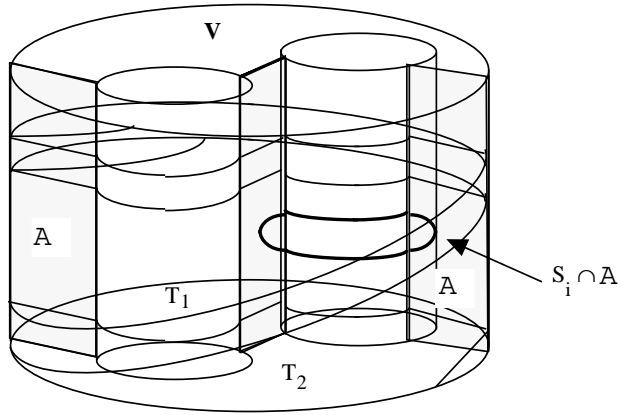

Figure 3.

**Step 2:** At the beginning of the sweep-out, no curve of intersection of $S_i \cap \partial T$ can be inessential in $\partial T$ and cut off a meridian of $B_i$ and at the end, no such curve can cut off a meridian of $A_i$. So again there is a subinterval where no inessential curve of $S_i \cap \partial T$ in $T$ can cut off a meridian in either $A_i$ or $B_i$. Restrict to such a subinterval, chosen to have the property that just before the beginning a meridian of $A_i$ is cut off and just after the end, a meridian of $B_i$ is cut off.

**Step 3:** It is still possible that during the sweep-out there will be curves in $S_i \cap \partial T$ that are inessential in both surfaces; those that lie entirely inside of $\mathcal{A}$ are easily removed, so we focus on those that slop across the annuli $\partial V \cap \partial T$. (See Figure 3.) Anytime during the sweep-out that there is such a curve in $S_i \cap \partial T$ (inessential in $T$ and $S_i$, and not lying entirely in $\mathcal{A}$) there will be an outermost inessential arc of $S_i \cap \mathcal{A}$, cutting off a disk that lies in $A_i$ or in $B_i$. This allows a $\partial$-compression of $S_i$ into $\partial V$ and thereby reveals a meridian disk in $A_i$ or $B_i$ lying in $V$.

If, at any stage during the sweep-out, there are simultaneously such $\partial$-compressions via $A_i$ and via $B_i$ then we could complete the proof of Theorem 3.8 immediately: If the $\partial$-compressions were into disjoint annuli in $\partial V - S_i$ it would contradict strong irreducibility, so in fact the annuli into which they $\partial$-compress must be adjacent in $\partial V$. More generally, even if the annuli are adjacent, the $\partial$-compressing disks themselves may be taken to be disjoint. So the boundary compressions can be done simultaneously, and the resulting surface still intersects $\partial V$ in essential curves. We can imagine $S_i$ isotoped so that both boundary compressing disks lie in a small collar of a component $\partial_0 V$ of $\partial V$. Put another way, we can isotope $S_i$ so that it intersects a collar of $\partial_0 V$ in an easily described way: There is a single horizontal or vertical tube attached to a collection of annuli which are either spanning annuli or $\partial$-parallel in the collar $\partial_0 V \times I$. In particular, this collar contains meridian disks of both $A_i$ and $B_i$. What remains of $S_i$ in $V$ when this collar is removed must then be incompressible, hence $\partial$-incompressible, by strong irreducibility. Hence the surface $S'$ obtained by compressing the single tube identified in the collar $\partial_0 \times I$ is either vertical or horizontal, and $S$ is obtained from $S'$ by attaching either (respectively) a horizontal or a vertical tube, as claimed by Theorem 3.8.

**Step 4:** Following the previous step, we may as well assume there are not simultaneously disjoint $\partial$-compressions into $A_i$ and $B_i$, so there is a subinterval in which there are no $\partial$ compressions at all, and this implies that in this subinterval, every curve of intersection in $S_i \cap \partial T$ is essential in $\partial T$, as was our claim.

We now want to understand how $S_i$ intersects each of the tori in $T$. This question is well-understood (see [Sc2], [MR]). There are four possibilities for each solid torus $T_j \in T$:

- $S_i$ intersects $T_j$ in meridian disks
- $S_i$ intersects $T_j$ in meridian disks and exactly one $\partial$-parallel annulus, parallel to a meridinal annulus in $\partial T_j$.
- $S_i$ intersects $T_j$ in a family of incompressible annuli in $T_j$
- $S_i$ intersects $T_j$ in a family of incompressible annuli, plus one other component obtained by tubing two incompressible annuli together or one annulus to itself, via a $\partial$-parallel tube.

If any $T_j$ contains a component of the second type, we can, by the previous argument, push a tube to the outside of $V$ and what remains inside will be incompressible, as required. If any $T_j$ contains a component of the fourth type, we are done by the same argument, unless the incompressible annuli are vertical, since if the annuli are vertical we do not know that they will $\partial$-compress into $\partial V$. Just as we eliminated inessential curves of intersection with $\partial T$ earlier, a further subinterval of the sweep-out can be found in which the third type does not arise, except perhaps when the annuli are vertical and so do not necessarily $\partial$-compress to $\partial V$. So the only remaining possibilities are the first, and also the third and fourth when the annuli are vertical.

Clearly meridian disks and vertical annuli cannot occur in neighboring solid tori, since their boundaries would intersect in some annulus of $\mathcal{A}$. It follows that either $S_i$ intersects each $T_j$ in meridian disks (i. e. $S_i$ intersects $V$ in a horizontal incompressible surface) or $S_i$ intersects each $T_j$ in vertical annuli, plus possibly somewhere a single horizontal tube. (More than one would contradict strong irreducibility.) This last would mean that $S_i$ intersects all of $V$ in vertical annuli and tori, with possibly one horizontal tube attached. □

**Definition 3.9.** *Suppose $V$ is a Seifert manifold. Then any manifold obtained by replacing a tubular neighborhood of a regular fiber by an irreducible, $\partial$-incompressible, acylyndrical, atoroidal manifold (with torus boundary) is called a* scrambling *of $V$.*

**Corollary 3.10.** *Suppose, in the conclusion of Theorem 3.8, there is a component of $S_i \cap V$ that is obtained from vertical surfaces by attaching a single horizontal tube. Then we can scramble $V \subset M$ so that the resulting manifold $M'$ has a Heegaard splitting of the same genus as $M$.*

**Proof:** The core of the horizontal tube projects to an arc in $P$ and that arc has preimage a vertical annulus $\alpha \subset V$. Suppose the tube compresses via a disk $D_A$ in $A_i$, say. Then the complement of the tube in $\alpha$ is a compressing disk $D_B$ for $B_i$ and $\partial D_B$ intersects $\partial D_A$ in two points. (See Figure 4.) A neighborhood of $\alpha$ in $V$ is a solid torus $T$ whose boundary intersects $S_i$ in four vertical curves and both $A_i \cap T$ and $B_i \cap T$ are genus two handlebodies. Viewing these genus two handlebodies as balls with two unknotted arcs removed, another view of $(A_i \cup_{S_i} B_i) \cap T$ is that $T$ is naturally homeomorphic to the complement of the unknot $u \subset S^3$, and, if we put the unknot in 2-bridge position, then $S_i$ is a 4-punctured equatorial sphere dividing $S^3 - u$ into two balls with two unknotted arcs in each. Replace $T$ with some atoroidal acylyndrical 2-bridge knot complement. This scrambles $V$ but merely changes the attaching map of $\partial_+ A_i$ and $\partial_+ B_i$ along $S_i \cap T$. □

## 4. First results

We will be presenting three extended arguments, each one of greater complexity than the preceding, but yielding more refined results. In this section we give the easiest of the three, whose conclusion implies immediately (and says much more than) that $M - \Theta$ has at most $3g - 3$ components, at most $g - 1$ of which are not Seifert.

**Definition 4.1.** *Suppose $\mathcal{A} \subset H$ is a properly imbedded essential collection of annuli in a compression body $H$. A complementary component $Y$ is* toral *if it is either a solid torus or is homeomorphic to torus $\times I$.*

*If the toral component $Y$ is a solid torus, define the* complexity *of $Y$ to be $c(Y) = |\partial_+ H \cap Y| - \epsilon$, where $\epsilon = 1$ if the annuli $\mathcal{A} \cap \partial Y$ are longitudes and $\epsilon = 0$ if they are not. In the former case, say that $Y$ intersects $\mathcal{A}$* longitudinally *or $Y$ is longitudinal. Otherwise say $\mathcal{A}$ is* twisted *in $Y$, or $Y$ is twisted. (See Figure 5.)*

*If $Y \neq H$ is homeomorphic to torus $\times I$ then define $c(Y) = |\partial_+ H \cap Y|$.*

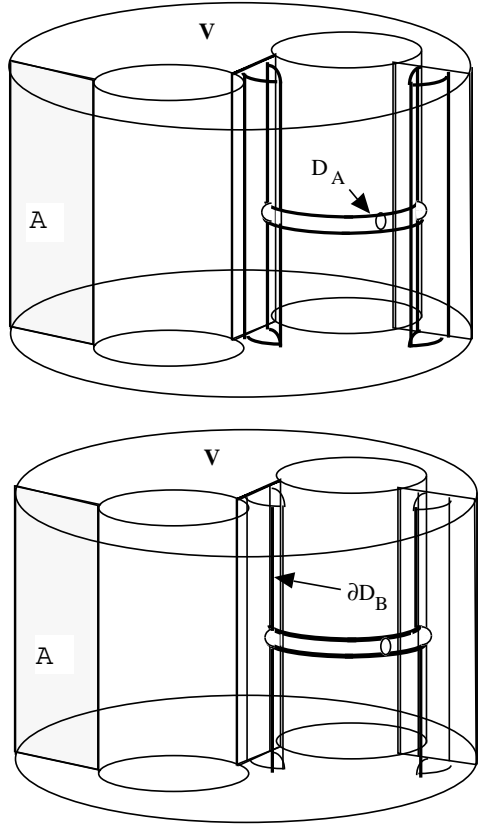

Figure 4.

In the degenerate case in which $Y = H$ is itself just torus $\times I$, and $\mathcal{A} = \emptyset$, say $c(Y) = 0$.

The complexity of the union of toral components is the sum of the complexities of each individual component.

Note that no essential annulus can have both boundary components in $\partial_- H$, so it follows that any solid torus component $Y \subset H$ intersects $\partial_- H$ in at most one component, an annulus. Such an intersection occurs if and only if exactly two components of $\mathcal{A} \cap \partial Y$ are spanning annuli. We then know more:

**Lemma 4.2.** *If a solid torus complementary component $Y$ of $\mathcal{A} \subset H$ intersects $\partial_- H$ in a single annulus, then $Y$ intersects $\mathcal{A}$ longitudinally.*

**Proof:** : Let $\alpha \subset \partial_- H$ denote the core of the annulus of intersection. If $\alpha$ is non-separating, consider a spanning annulus whose end in $\partial_- H$ is a simple closed curve which intersects $\alpha$ in a single point. By standard innermost disk, outermost arc arguments, we can choose the annulus so that it intersects $Y$ in a single disk $D$ and $\partial D$ crosses $\alpha$ exactly once. This guarantees that $D$ is a meridian and so $\alpha$ is a longitude of $Y$.

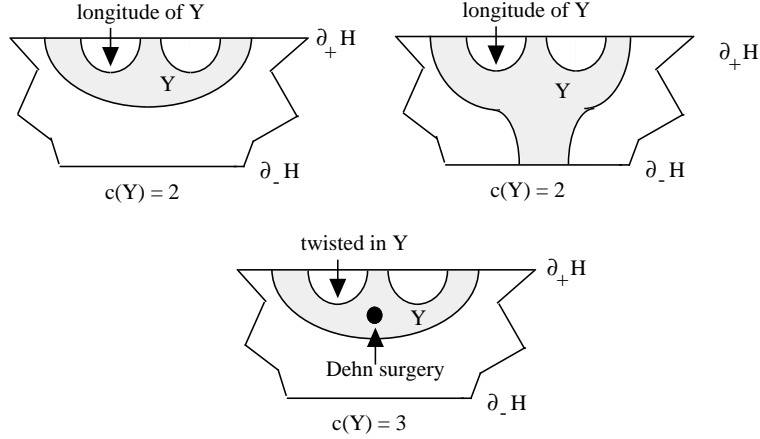

FIGURE 5.

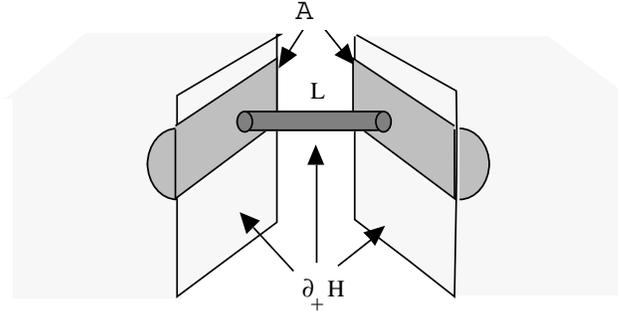

FIGURE 6.

If $\alpha$ is separating, the argument is only slightly more complex: The same sort of innermost disk, outermost arc argument shows that one can find spanning annuli $A_1$ and $A_2$ so that each curve $\alpha_i = A_i \cap \partial_- H$ is parallel to $\alpha$ but the two curves $\alpha_1, \alpha_2$ lie on opposite sides of $\alpha$ in $\partial_- H$, and so that $A_i \cap Y = \emptyset$. Then choose $\gamma$ to be a spanning arc for the annulus between the $\alpha_i$ in $\partial_- H$ and construct a "spanning square" $\Sigma \subset H$ so that $\Sigma \cap \partial_- H = \gamma$, $\Sigma \cap A_i$ is a spanning arc of $A_i$ and the rest of $\partial \Sigma$ lies in $\partial_+ H$. Then an innermost disk, outermost arc argument on $\Sigma \cap \mathcal{A}$ shows that $\Sigma$ can be chosen so that $\Sigma \cap Y$ is a meridian disk $D$ as above. □

**Definition 4.3.** *Suppose $\mathcal{A} \subset H$ is a properly imbedded essential collection of annuli in a compression body $H$. A complementary component $L$ is a* basic block *if it is a genus two handlebody in which at least two components of $\mathcal{A} \cap L$ are longitudes, separated by a meridian disk of $L$.*

*A complementary component is a* spanning product *if, for some surface $E$ it is of the form*
$$E \times (I; 0, 1) \subset (H; \partial_- H, \partial_+ H).$$

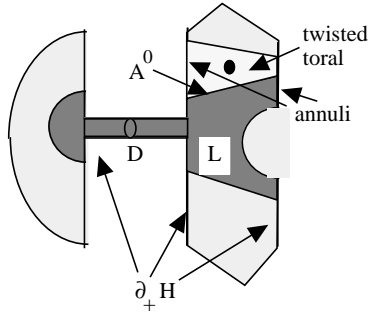

Figure 7.

**Lemma 4.4.** *Suppose $H$ is a compression body, and $\mathcal{A}$ is a properly imbedded essential collection of annuli in $H$. Suppose no two adjacent components of $H - \mathcal{A}$ are toral, and there are $n$ non-toral components of $H - \mathcal{A}$ that are disjoint from $\partial_- H$. Then*

$$n \leq J(H)/2.$$

*Moreover, if $Y$ is the union of toral components of $H - \mathcal{A}$ and $\alpha$ is the number of non-spanning annuli which are not adjacent to toral components (on either side), then*

$$c(Y) + \alpha \leq J(H).$$

*Finally, suppose in fact $c(Y) + \alpha = J(H)$. Then*

- $n = J(H)/2$.
- *Each non-toral component is either a basic block or a spanning product.*
- *If more than one annulus of $\mathcal{A}$ lies on the same side of a separating meridian disk $D$ for a basic block $L$ then one of those annuli $A^0$ has either one non-toral component of $H - \mathcal{A}$ on each of its sides in $H$ or one annulus component of $\partial_+ H - \mathcal{A}$ on each of its sides. (See Figure 7.)*

**Proof:** We will proceed by induction on the pair $(J(H), c(Y))$.

When $J(H) = 0$ then $H$ is a product. In a product, the only essential annuli are spanning annuli. Then $J(H) = c(Y) = \alpha = 0$, every complementary component is a spanning product, and the lemma is true in this case.

Note that any toral component $Y^0$ of $H - \mathcal{A}$ that has $c(Y^0) = 0$ is just a product neighborhood of a spanning annulus. The two annuli in $\mathcal{A} \cap Y^0$ could be removed with no effect on the argument. So we may as well assume that any toral component has positive complexity.

Suppose that some component $Y^0$ of $H - \mathcal{A}$ is a solid torus component of $H - \mathcal{A}$ (with positive complexity) that intersects $\mathcal{A}$ longitudinally. If there is a spanning annulus in $\mathcal{A} \cap \partial Y^0$ delete it from $\mathcal{A}$; otherwise delete a non-spanning annulus in $\mathcal{A} \cap \partial Y^0$ from $\mathcal{A}$. This attaches $Y^0$ to a non-torus component of $H - \mathcal{A}$, thereby lowering $c(Y)$ by $c(Y^0)$. On the other hand, every non-spanning component of $\mathcal{A} \cap \partial Y^0$ (except, perhaps, the one removed) now contributes to $\alpha$. So $c(Y)$ is lowered, but $\alpha + c(Y)$ remains the same. (See Figure 8.) The operation leaves both $J(H)$ and the number of non-toral components not intersecting $\partial_- H$ the same. It cannot create a new basic block or spanning product. On any basic block the operation may increase but cannot decrease the number of annuli intersecting the basic block. The operation can decrease but cannot increase the number of annuli in $\partial H - \mathcal{A}$. So the result follows by induction.

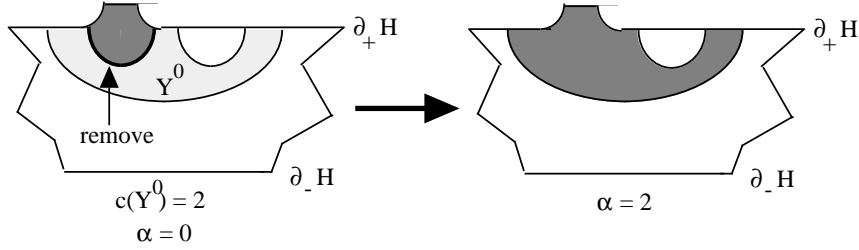

FIGURE 8.

There remains only the case in which $J(H) > 0$ and no solid torus component of $H - \mathcal{A}$ intersects $\mathcal{A}$ longitudinally. In particular, no spanning annulus is adjacent to any torus component. Unless $\mathcal{A}$ is $\partial$-compressible, $\mathcal{A}$ would consist entirely of vertical spanning annuli, so $\alpha = c(Y) = 0$, and the result follows immediately. So we will suppose that $\mathcal{A}$ is $\partial$-compressible. Let $A \in \mathcal{A}$ be a component of $\mathcal{A}$ on which such a $\partial$-compression can be done, and let $D$ be the disk that results from $A$ by the $\partial$-compression.

First suppose there is such a choice of $\partial$-compression so that $D$ is separating in H. Cut $H$ open along $D$ to get two components: $H_1$ which comes from the side of $A$ that contains $D$, and $H_2$ which comes from the opposite side. Let $\mathcal{A}_1 = \mathcal{A} \cap H_1$ and $\mathcal{A}_2 = (\mathcal{A} - A) \cap H_2$. Note that no component of $\mathcal{A}_2$ can be an inessential annulus, for if any were, it would have been adjacent to a solid torus component of $H - \mathcal{A}$ in which $A$ is longitudinal. So, by induction, the lemma is true for $\mathcal{A}_2 \subset H_2$ (or, possibly, $H_2$ is simply a solid torus on which $A$ is not longitudinal). On the other hand, we cannot immediately assume that the inductive hypothesis applies to $H_1$, for two things could go wrong. First, it's quite possible that exactly one component $A_1$ of $\mathcal{A}_1$ is inessential in $H_1$, exactly when $A \cup A_1$ together cut off a genus two solid handlebody component of $H - \mathcal{A}$ (containing $D$) in which both annuli $A$ and $A_1$ are longitudinal, i. e. a basic block $L$ with $|L \cap \mathcal{A}| = 2$. Another problem could be that the component $Y_0$ adjacent to $D$ in $H_1$ is a solid torus, and also $Y_0$ is adjacent to other toroidal components, in violation of the inductive hypothesis. Note that in this case $Y_0$ must be longitudinal and adjacent via only one annulus (which we again call $A_1$) to a toroidal component, since any annulus in $H$ is $\partial$-compressible and all solid torus components of $H - \mathcal{A}$ are twisted. The component of $H - \mathcal{A}$ containing $D$ is a basic block $L$, but this time $|L \cap \mathcal{A}| \geq 3$, and furthermore $A_1$ intersects an annulus component of $\partial_+ H - \mathcal{A}$ on both sides. In either case, removing $A_1$ from $\mathcal{A}_1$ (call the result $\mathcal{A}'_1$) restores the inductive hypothesis to $\mathcal{A}'_1 \subset H_1$.

If either failure of the inductive hypothesis arises (so that the component containing $D$ is a basic block $L$, and an incident annulus $A_1$ is either inessential in $H_1$ or adjacent to a toroidal component $Y_1$ of $H - \mathcal{A}$) then, back in $H$, removing $A_1$ from $\mathcal{A}$ would either decrease $\alpha$ by one, or decrease $c(Y)$ by $c(Y_1)$ and increase $\alpha$ by $c(Y_1) - 1$. Similarly, if the annulus $A$ is adjacent only to non-toral components, then removing it would decrease $\alpha$ by one; if it is adjacent to a solid torus $Y_2$, removing $A$ would lower $c(Y)$ by $c(Y_2)$ and raise $\alpha$ by $c(Y_2) - 1$. In either case, removing the annulus would decrease $c(Y) + \alpha$ for $\mathcal{A} \subset H$ by exactly one. So if and only if $A_1$ exists (i. e. the inductive hypothesis fails) we get that $c(Y) + \alpha$ drops by two in going from $\mathcal{A} \subset H$ to $\mathcal{A}'_1 \cup \mathcal{A}_2$ in $H_1 \cup H_2$ (which is necessary if the second inequality is an equality). So if the second inequality is an equality, $A_1$ exists.

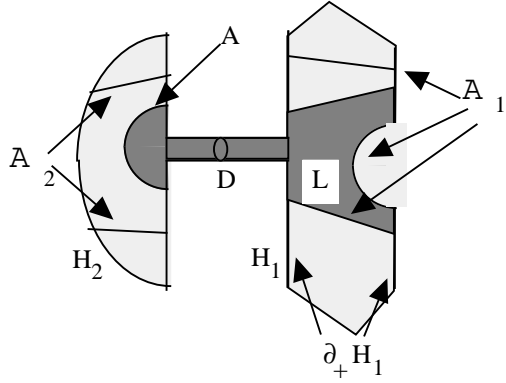

FIGURE 9.

Note that $J(H) = J(H_1) + J(H_2) + 2$ and $H - \mathcal{A}$ has at most one more non-toral component not intersecting $\partial_- H$ than the sum of the number of such components in $H_1$ and $H_2$. It has exactly one more unless the component adjacent to $D$ in $H_1$ is itself not toral, with the consequence that $A_1$ is not defined. So the proof now follows by induction: If the second inequality is an equality, then $A_1$ exists, the component $L$ containing $D$ is a basic block satisfying the required conditions (as we have seen), and also the second inequality must remain an equality for $\mathcal{A}'_1 \cup \mathcal{A}_2 \subset H_1 \cup H_2$. There we apply the inductive hypothesis to show that all non-toral components other than the one containing $D$ are spanning products or basic blocks satisfying the required conditions, and the first inequality for $\mathcal{A}'_1 \cup \mathcal{A}_2 \subset H_1 \cup H_2$ is an equality. We have just shown that if $A_1$ exists, then this equality implies the first equality for $\mathcal{A} \subset H$.

Finally, suppose that such a $\partial$-compression results in a *non-separating* disk $D$, so $D$ is incident to a non-separating annulus $A \subset \mathcal{A}$. Cut $H$ open along $A$ and, along the copy of $A$ that is still incident to $D$, attach a twisted solid torus to get a new handlebody $H'$ (since the attaching annulus is longitudinal after the cut) and a family of annuli $\mathcal{A}' \subset H'$ with the same number of elements as $\mathcal{A}$. Indeed we can think of just replacing $A$ with a copy of $A$. It's easy to see that neither operation affects $J(H)$ or the number or type of non-toral components. The cut along $A$ reduces one of $c(Y)$ or $\alpha$ by exactly one, depending on whether or not on the other side of $A$ from $D$ is a (twisted) toral component or a non-toral component. On the other hand, attaching the twisted torus raises $c(Y)$ by exactly one. So the result for $\mathcal{A}' \subset H'$ would imply the result for $\mathcal{A} \subset H$. But a proof for $\mathcal{A}' \subset H'$ is provided by the previous case, since $A$ is separating in $H'$. □

**Definition 4.5.** *Let $V$ be a connected orientable Seifert manifold with non-empty boundary, $f$ singular fibers, and orbit space the surface $P$. Define the* vertical index *of $V$ to be*
$$I_v(V) = f - \chi(P).$$
*Define the augmented vertical index as*

$$I_v^+(V) = \begin{cases} 1 + I_v(V) = 1 - \chi(P) & \text{if } V \text{ has no exceptional fibers.} \\ 1/2 + I_v(V) = 3/2 - \chi(P) & \text{if } V \text{ has one exceptional fiber} \\ I_v(V) = f - \chi(P) & \text{if } V \text{ has } f \geq 2 \text{ exceptional fibers} \end{cases}$$

*Define the epsilon vertical index as*

$$I_v^\epsilon(V) = \begin{cases} I_v^+(V) & \text{when } |\partial V| = 1 \\ I_v^+(V) - 1/2 = I_v(V) + 1/2 = 1/2 - \chi(P) & \text{if } |\partial V| = 2 \text{ and } f = 0 \\ I_v(V) = f - \chi(P) & \text{otherwise} \end{cases}$$

*If $V$ is not connected, the (resp. augmented, epsilon) vertical index of $V$ is defined to be the sum of the (resp. augmented, epsilon) vertical indices of its components.*

Note:
- If $I_v^+(V) \leq 1$ or $I_v^\epsilon(V) \leq 1$ then $V$ is atoroidal. Indeed, the only toroidal manifolds for which either index is $\leq 2$ are those for which $f = 1$ and $P$ is the Möbius band, or $f = 0$ and $P$ is the once-punctured Möbius band.
- A scrambling of $V$ effectively removes from the Seifert piece a vertical solid torus (and creates a non-Seifert piece). So scrambling deletes a disk from $P$ and so the process will increase both $I_v(V)$ and $I_v^+(V)$ by exactly one.

The theorems below will put bounds on the various complexities, based on the genus of the Heegaard splitting. It would make the definitions considerably easier, the aesthetics better, and the theorems stronger, if the fractions in the definitions of augmented and epsilon index above could be raised to the nearest integer. We know of no counterexamples to this hope. The difference between the types of complexity above are small, and would not be worth making, except that, following Kobayashi's agenda, we would like to be able to get a bound on the number of toroidal components and a few toroidal manifolds have unaugmented vertical index 1: a circle bundle over the once punctured torus or Klein bottle, or over the punctured Möbius band, or a Seifert bundle over the Möbius band with one exceptional fiber. For this reason we need to consider augmentation.

**Definition 4.6.** *Suppose $M$ is a closed irreducible 3-manifold, $V$ is a Seifert component from the torus decomposition of $M$, and $M$ has the (generalized) Heegaard splitting*
$$M = (A_1 \cup_{S_1} B_1) \cup_{F_1} (A_2 \cup_{S_2} B_2) \cup_{F_2} ... \cup_{F_{m-1}} (A_m \cup_{S_m} B_m).$$

Let $S = \cup S_i$ and $F = \cup F_i$. Then $V$ is aligned *with respect to the Heegaard splitting* if $S \cup F$ intersects $\partial V$ only in vertical fibers.

Note that if the canonical tori $\Theta$ have been put in preferred position with respect to the splitting, then a Seifert component $V$ is aligned if *any* boundary component of $V$ intersects any surface $F_i, S_i$ in a vertical fiber. If $F \cup S$ is disjoint from a boundary component $T$ of $V$ then, since $\Theta$ is incompressible, it must be that $T$ is parallel to a component of $F$. If this happens for any component of $\partial V$ then this, too, implies that $V$ is aligned and we say that that component of $\partial V$ is *strongly aligned*.

**Theorem 4.7.** *Suppose $M$ is a closed orientable $3$-manifold that has a genus $g$ irreducible Heegaard splitting. Let $\Theta$ be the collection of canonical tori for $M$, put in preferred position with respect to a strongly irreducible generalized Heegaard splitting for $M$, and suppose no component of $\Theta$ is strongly aligned. Denote by $V$ the union of aligned Seifert components and let $p$ be the number of components of $M - \Theta$ that are not aligned Seifert components. Then*

- $p \leq g - 1$ and
- $I_v^+(V) \leq 2g - 2$.

**Proof:** Suppose that
$$M = (A_1 \cup_{S_1} B_1) \cup_{F_1} (A_2 \cup_{S_2} B_2) \cup_{F_2} ... \cup_{F_{m-1}} (A_m \cup_{S_m} B_m)$$
is a strongly irreducible splitting, of genus $\leq g$, with respect to which the canonical tori $\Theta$ can be put in preferred position with no components of $\Theta$ strongly aligned. There may be some of the $S_i$ which intersect the aligned Seifert manifolds in vertical surfaces to which a horizontal tube has been attached, as described in 3.8. The proof will be by induction on the number of such $S_i$. The inductive step is easy, so we present it before examining the case when there are no such horizontal tubes.

If there is a component of $S_i \cap V$ which contains a horizontal tube, apply Corollary 3.10 to scramble $V$, replacing $M$ with a manifold $M'$ of the same genus but one for which the number of horizontal tubes in aligned Seifert pieces has been reduced by one, the number of components which are not aligned Seifert pieces is increased by one, as is $I_v^+(V)$. The theorem for $M'$ then implies the theorem for $M$.

It remains to prove the theorem under the assumption that there are no horizontal tubes of any $S_i$ in $V$, so all components of $(F \cup S) \cap V$ are vertical surfaces in $V$, either incompressible tori or essential annuli. Hence each component of each $A_i \cap V$ and $B_i \cap V$ is toral.

We first count the number of components of $M - \Theta$ that cannot be viewed as aligned Seifert manifolds. Suppose $N$ is such a component. Then no component of intersection of $N$ with any $A_i$ or $B_i$ is toral, for if one were then $N$ would contain an essential annulus and it could not be vertical; this would imply that $N$ is a Seifert piece (since it contains an essential annulus) and the annulus would be horizontal (by Theorem 3.5). Hence $N$ is an $I$-bundle construct with $I$-base the Möbius band or the annulus. In either case, it could be fibered differently and so be viewed as aligned.

The upshot is that, to get a bound on the number of components that are not aligned Seifert components, we need to count exactly those that do not intersect any compression body in any toral component. Moreover no two adjacent components of $A_i - \Theta$ (or $B_i - \Theta$) can be toral, since if they were, then the component of $\Theta$ between them would be incident to aligned Seifert pieces on both sides, and so should have been eliminated from $\Theta$. It follows then from 4.4 that the number of components of $A_i - \Theta$ that are not part of an aligned Seifert piece and which don't intersect $\partial_- A_i$ is at most $J(A_i)/2$. So this is an upper bound to the number of components of $M - \Theta$ that simultaneously are not aligned Seifert pieces, intersect $A_i$, and are disjoint from all $A_j, B_j, j \leq i$. On the other hand, any component of $M - \Theta$ which is not an aligned Seifert component must intersect some $A_i$ and some $B_i$ and do so in non-toral pieces. So starting at $A_1$ and summing across each part of the generalized Heegaard splitting, we conclude that there are at most $\sum_{i=1}^{m} J(A_i)/2 = g - 1$ components which are not aligned Seifert components. This verifies the first claim.

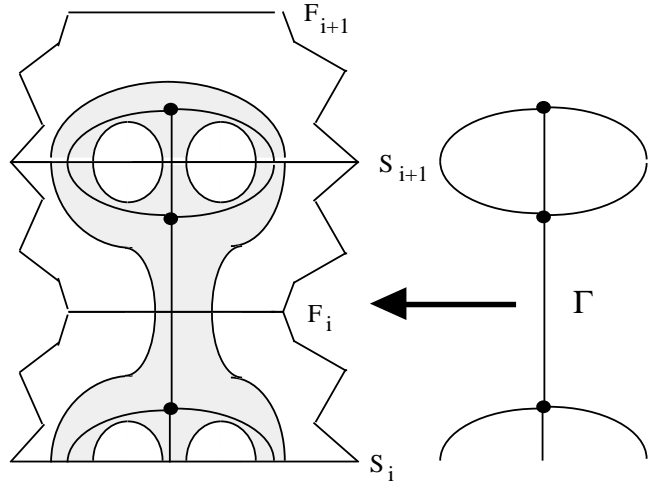

FIGURE 10.

Now we count the index of the aligned Seifert piece $V$. Let $P$ be the base surface. We construct a graph $\Gamma$ which is a deformation retract for $P$, noting exceptional fibers as we go (See Figure 10): Include a vertex for every solid torus component of $A_i \cap V$ and $B_i \cap V$. (Call the vertex *longitudinal* if the annuli in which the torus intersects $S_i$ are longitudinal.) Include a loop for every *torus × I* component of the intersection. Include an edge in $\Gamma$ for every annulus in $F_i \cap V$ and $S_i \cap V$, with ends of the edge on the vertices corresponding to the solid torus (or *torus × I*) component on either side of the annulus. The exceptional fibers correspond precisely to the $f$ torus components whose intersection with $S_i$ is not longitudinal, i. e. the $f$ non-longitudinal tori. If $v$ is the number of vertices in $\Gamma$ and $e$ is the number of edges (ignoring loops in both counts), then $I_v(V) = e - (v - f)$.

Now $e$ is half the sum of the valences of the vertices, so a local formula for $I_v(V)$ can be gotten as follows: For each vertex $u$ in $\Gamma$ not the base of a loop, note that $c(\tau) = valence(u)$ if the corresponding solid torus $\tau$ (in some $A_i \cap V$ or $B_i \cap V$) is not longitudinal; $c(\tau) = valence(u) - 1$ when $\tau$ is longitudinal and disjoint from the $F_i$; and $c(\tau) = valence(u) - 2$ when $\tau$ intersects the $F_i$. For each loop representing a component $\alpha \cong torus \times I$, $c(\alpha)$ denotes the number of ends at the base of the loop not belonging to the loop itself. Abusing notation a bit, for $u$ a vertex (not the base of a loop) in $\Gamma$ representing $\tau$, let $c(u) = c(\tau)$ and for $\alpha \subset \Gamma$ a loop, let $c(\alpha)$ denote the complexity of the corresponding component $\cong (torus \times I)$. Then

$$2I_v(V) = \sum_{u \in \Gamma} valence(u) + \sum_{\alpha \subset \Gamma} c(\alpha) - 2(v - f)) = \sum_{u \in \Gamma} c(u) + \sum_{\alpha \subset \Gamma} c(\alpha) - \Lambda,$$

where $\Lambda$ is the total number of longitudinal components occuring in all $A_i \cap V$ or $B_i \cap V$ that are disjoint from $F$.

To summarize, if we let $a_i$ (resp. $b_i$) represent the sum of the complexities of all the toral components of $A_i \cap V$ (resp. $B_i \cap V$) then

$$(1) \qquad 2I_v(V) = \sum_{i=1}^{m}(a_i + b_i) - \Lambda.$$

No boundary component of $V$ is totally aligned, by assumption. This means that for any component $V^0$ of $V$ and the lowest $i$ for which $A_i \cap V^0$ is non empty, in fact $A_i \cap V^0$ must be a solid torus disjoint from $\partial_- A_i$. Similarly for the largest $i$ for which $B_i \cap V^0$ is non empty. So the component $V^0$ contributes at least 2 to $\Lambda$ if it has no exceptional fibers and contributes at least 1 to $\Lambda$ if it has only one exceptional fiber. The upshot is that we have the inequality:

$$2I_v^+(V) \leq \sum_{i=1}^{m}(a_i + b_i) = \sum_{i=1}^{m} a_i + \sum_{i=1}^{m} b_i$$

On the other hand, Lemma 4.4 above shows that

$$\sum_{i=1}^{m} a_i + \sum_{i=1}^{m} b_i \leq \sum_{i=1}^{m} J(A_i) + \sum_{i=1}^{m} J(B_i) = 4g - 4.$$

Hence we get the second inequality

$$I_v^+(V) \leq 2g - 2.$$

□

**Corollary 4.8.** *If, in the conclusion of Theorem 4.7, we have $I_v^+(V) = 2g - 2$, then each of $M - \Theta$ is either an aligned Seifert manifold or a 2-bridge link exterior.*

**Proof:** In order to achieve equality, all inequalities in the proof must be equalities. This implies that when Lemma 4.4 is applied to each compression body $A_i$ (or $B_i$) we have $\alpha = \emptyset$, and each non-toral component is a basic block or a spanning product. Since $\alpha = 0$ it follows immediately from the proof of Lemma 4.4 (in which longitudinal solid tori are traded in immediately for annuli in $\alpha$) that the only longitudinal solid tori are spanning products. In particular, any annulus of $A^0 \in \mathcal{A} = \Theta \cap A_i$ that is incident to a basic block is incident on the other side to a twisted solid torus. In particular, since $A^0$ must $\partial$-compress, it will $\partial$-compress through the basic block. This boundary compressing disk $D$ can be chosen to avoid all other components of $A_i - \mathcal{A}$, since of those that it might intersect, that one adjacent to $L$ cannot be non-toral (since $\alpha = 0$), nor a longitudinal torus (we have eliminated that possibility too), or a twisted torus (since $D$ intersects one of its boundary annuli in exactly one arc). Hence each annulus in each basic block of $A_i$ $\partial$-compress within the basic block to $\partial_+ A_i$.

In particular we conclude, (or could also from the last claim of 4.4) that each basic block $L$ in $A_i$ has exactly one component of $\mathcal{A}$ on each side of a separating meridian disk. Since $L$ was chosen so that one of the annuli $\partial$-compresses through $L$, $\partial L \cap \partial H$ is a 4-punctured sphere.

Any $\partial$-compression in a basic block in $A_i$ must be into the same component of $\partial_+ A_i - \Theta$ as a $\partial$-compression into $\partial_+ B_i - \Theta$, since $A_i \cup_{S_i} B_i$ is strongly irreducible. This means that there can be at most one basic block in $A_i$ and one in $B_i$ and they are glued together along the 4-punctured sphere on which they coincide along $S_i$. All other components must be toral. □

**Corollary 4.9.** *Suppose $M$ is a closed orientable 3-manifold that has a genus $g$ irreducible Heegaard splitting. Let $\Theta$ and $V$ be as in Theorem 4.7 but allow the possibility that some components of $\partial V$ are strongly aligned. Then $I_v^\epsilon(V) \leq 2g - 2$.*

**Proof:** Suppose a component $T$ of $\partial V$ is parallel to a component of $F$. The proof is by induction on the number $t$ of strongly aligned components of $\partial V$. Since $I_V(V) \leq I_v^\epsilon(V) \leq I_v^+(V)$, the case $t = 0$ is proven in Theorem 4.7. So suppose some component $T \in \Theta$ of $\partial V$ is parallel to a component of $F$. Let $V^0$ denote (one of the) Seifert components on which $T$ lies.

Construct a new manifold $M'$ from $M$ by cutting $M$ open along $T$ and gluing on twisted solid tori to each copy of $T$. These attaching maps can be chosen so that if, on the other side of $T$ there is an acylindrical atoroidal piece of $M - \Theta$, it will remain so in $M' - \Theta'$. The effect on $V$ is to replace one or two boundary components by essential fibers. This has no effect on $I_v V$ but may cause a loss of augmentation in $V^0$. Hence we have immediately that $I_v(V) \leq 2g - 2$. What we need to show is that this loss can be limited if $|\partial V^0| \leq 2$.

Suppose $\partial V^0$ has only one boundary component $\partial V^0$ and it is strongly aligned. The manifold $M'$ then consists of two components, $V_+$ (of some genus $g_v$) which contains the Dehn filled $V^0$, and $M_-$ of some genus $g_-$ with $(g-1) = (g_- - 1) + (g_v - 1)$. Since $V^0$ is entirely Seifert, some compressible component of some $S_i$ lay entirely in $V^0$. This means, via 3.10, that $V^0$ could have been scrambled without affecting genus $M$ but raising $I_v^+(V_+)$ by 1. This would counterbalance any loss of augmentation; in fact we must have had $I_v^\epsilon(V) = I_v^\epsilon(V - V^0) + I_v^+(V^0) \leq I_v^\epsilon(V - V^0) + (I_v(V_+) + 1) \leq 2(g_- - 1) + 2(g_- - 1) = 2g - 2$.

If $\partial V^0$ has two boundary components and they are both totally aligned, the same induction argument shows that $I_v^\epsilon(V) \leq 2g - 2$, so we need only consider the case in which exactly one boundary component is totally aligned. If $V^0$ has exceptional fibers then, since $I_v^\epsilon(V^0) = I_v(V^0)$ so there is not augmentation and the proof follows immediately by induction from Dehn filling at the strongly aligned component, as above. If $V^0$ has no exceptional fibers, then after the Dehn filling, the augmentation drops from 1 to 1/2, but also $I_v^\epsilon(V^0) = I_v^+(V^0) - 1/2$, so again the argument follows by induction. □

## 5. Non-aligned Seifert pieces

In this section we are able to refine further so that we also have information about the Seifert components that are not aligned. The argument is, in outline, parallel to that of the previous section.

There is this technical refinement of 4.4.

**Lemma 5.1.** *Suppose $H$ is a compression body, $\mathcal{A}$ is a properly imbedded essential collection of annuli in $H$, and no two adjacent components of $H - \mathcal{A}$ are toral. Suppose $(E, \partial E) \subset (H, \partial_+ H)$ is a compact properly imbedded subsurface, with $s$ components and $d$ boundary components, and, among the components of $H - \mathcal{A}$, is a collection homeomorphic to an $I$-bundle $\xi$ over $E$, with $\dot{\xi} \subset \partial_+ H$ and $\xi|\partial \subset \mathcal{A}$. Then there are at most*

$$J(H)/2 + \chi(E) + s$$

*non-toral components of $H - \mathcal{A}$ that don't intersect $\partial_- H$. Moreover, if $Y$ is the union of the toral components of $H - \mathcal{A}$, then*

$$c(Y) \leq J(H) + 2\chi(E) + d.$$

**Proof:** Let $H'$ be the compression body (or bodies) obtained by deleting $\xi$ from $H$, and $\mathcal{A}' = \mathcal{A} - (\xi|\partial)$. Remove any further components of $\mathcal{A}$ that are $\partial$-parallel in $H'$. This latter move has no effect on the non-toral components of $H'$. Then $J(H') = J(H) + \chi(\dot{\xi}) = J(H) + 2\chi(E)$ and, since $\xi$ contains at most $s$ non-toral components, the assertion that there are at most $J(H)/2 + \chi(E) + s$ non-toral components of $H - \mathcal{A}$ that don't intersect $\partial_- H$ then follows from Lemma 4.4.

For the second result, first note that we may as well assume no component of $E$ is an annulus or Möbius band. For we could regard such a component as toral and the theorem in that case would suffice, for in no longer regarding the component as toral, $c(Y)$ goes down by 1, $s$ goes up by one, $J(H) + 2\chi(E)$ remains the same, and $d$ goes up by two.

Achieving the inequality $c(Y) \leq J(H) + 2\chi(E) + d$ is relatively simple. For this, let $H'$ be the compression body obtained by removing $\xi$ plus all toral components whose boundary intersects $\mathcal{A}$ exactly in $\xi|\partial$. Further remove from $\mathcal{A}'$ any annulus which is $\partial$-parallel in $H'$; each corresponds to a longitudinal toral component of $H - \mathcal{A}$ whose boundary contains only a single annulus in $\mathcal{A}'$. The effect is to reduce $c(Y)$ by at most $d$ and to set $J(H') = J(H) + 2\chi(E)$. ($c(Y)$ is not necessarily reduced exactly by $d$, because of longitudinal tori whose boundaries are contained in $\xi|\partial$ and components of $\xi|\partial$ which are not adjacent to toral components.) Then Lemma 4.4 applied to the compression body (or bodies) $H'$ completes the proof. □

**Definition 5.2.** *Let $V$ be an orientable Seifert manifold with non-empty boundary. Suppose $(E, \partial E) \subset (V, \partial V)$ is an incompressible possibly 1-sided surface in $V$ such that $V$ is an $I$-bundle construct over $E$). For $\rho$ the set of slopes in $\partial V$ determined by $\partial E$, call $I_\rho(V) = |V| - 2\chi(E) - |\partial E|$ the* horizontal index *of $V$ (with slope $\rho$).*

**Theorem 5.3.** *Suppose $M$ is a closed orientable 3-manifold that has a genus $g$ irreducible Heegaard splitting. Let $\Theta$ be the collection of canonical tori for $M$, put in preferred position with respect to a strongly irreducible generalized Heegaard splitting for $M$. Suppose $M - \Theta$ has $p$ non-Seifert components. Divide the Seifert components into two classes: the non-aligned or* horizontal *components $U_1, ..., U_h$ whose union we denote $U$, and the aligned Seifert components whose union we denote $V$. We know that $U$ can be written as as an $I$-bundle construct over a surface $E$ (each component of which has negative Euler characteristic). Denote the slope of $\partial E$ in $\partial U$ by $\rho$. Then*

- $p \leq g - 1 + \chi(E)$
- $I_\rho(U) + I_v^\epsilon(V) \leq 2g - 2 + h$.

**Proof:** Suppose that

$$M = (A_1 \cup_{S_1} B_1) \cup_{F_1} (A_2 \cup_{S_2} B_2) \cup_{F_2} ... \cup_{F_{m-1}} (A_m \cup_{S_m} B_m)$$

is the strongly irreducible splitting with respect to which $\Theta$ has been put in preferred position. By an inductive argument as in the proof of Theorem 4.7, we may as well assume that each component of each $A_i \cap V$ and $B_i \cap V$ is toral. Following the argument of Corollary 4.9, we may as well aim to prove the stronger inequality $I_\rho(U) + I_v^+(V) \leq 2g - 2 + h$ but also assume that no component of $\partial V$ is strongly aligned.

We first do a refined count of the number of components of $M - \Theta$ that are not aligned Seifert manifolds. As in the proof of Theorem 4.7, we have to count exactly those components that do not intersect any compression body torally.

It follows from 3.8 that, with possibly the exception of one component, the horizontal Seifert part intersects each $A_i$ in an $I$-bundle $\xi$. The possible exception, called a *tubed* component, starts with the form $\xi$, but then either a vertical tube along an $I$-fiber is deleted, or a $1 - handle$ is attached to $\dot{\xi}$. Ignore for the moment the possibility of a tubed component, and note then that the components of intersection of the horizontal Seifert pieces with each $A_i$ can be divided into two types: those that are spanning (i. e. are homeomorphic to $surface \times (I : 0, 1) \subset (A_i; \partial_- A_i, \partial_+ A_i)$) and those that are non-spanning and so are homeomorphic to $(\xi, \dot{\xi}) \subset (A_i, \partial_+ A_i)$. The union of the latter can be written in the form $(\xi_i, \dot{\xi}) \subset (A_i, \partial_+ A_i)$, where $\xi_i$ is an $I$-bundle over a surface $E_i$, just as in the hypothesis of Lemma 5.1.

For each component $U_j$ of $U$, choose an $i$ (there must be one, since $M$ is closed) for which at least one component of $U_j \cap A_i$ is an $I$-bundle as in Lemma 5.1. Similarly choose an $i'$ for which at least one component of $U_j \cap B_{i'}$ is an $I$-bundle as in Lemma 5.1. Then for each $1 \leq i \leq m$ let $h_i$ denote the number of components of $U$ that are assigned to $i$, $E_i$ denote the union of $I$-bases of those components of $U$ assigned to $i$ and $d_i$ represent the number of boundary components of $E_i$. Similarly let $h_{i'}$ denote the number of components of $U$ that are assigned to $i'$, $E_{i'}$ the union of the $I$-bases of those components and $d_{i'}$ the number of boundary components of $E_{i'}$. We repeat, each component of $U$ has its $I$-base surface assigned to exactly one of the $E_i, 1 \leq i \leq m$ and one of the $E_{i'}, 1 \leq i' \leq m$. With this set-up we are in a position to use the fact (from 5.1) that there are at most $J(A_i)/2 + \chi(E_i) + h_i$ non-toral components of $A_i - \Theta$ that don't intersect $\partial_- A_i$. Note that this means there are at most $J(A_i)/2 + \chi(E_i)$ components of $A_i - \Theta$ which are not parts of Seifert pieces, aligned or unaligned. We apply this fact, as we did the simpler inequality from Lemma 4.4 in the proof of 4.7 and conclude that there are at most

$$\sum_{i=1}^{m}(J(A_i)/2 + \chi(E_i)) = g - 1 + \chi(E)$$

non-Seifert components, the first inequality.

If there is a tubed component of $U \cap A_i$, it only improves the situation. Suppose first that there is a component of intersection of $A_i$ with a horizontal Seifert piece and it has the form $\xi \cup (1-handle)$. Then $\partial$-reducing the 1-handle returns us to the situation analyzed above, and reduces $J(A_i)$ by two. Lemma 5.1 after the $\partial$-reduction implies the required inequality before the $\partial$-reduction. Similarly, if a component is of the form $\xi - (fiber)$ then removing the component from $A_i$ reduces $J(A_i) + 2\chi(E_i)$ by two and both $h_i$ and the number of non-toral components by one. So again Lemma 5.1, applied after the removal of the component, implies the required inequality before the removal.

The second inequality is obtained in exactly the same way, using the inequality $c(Y) \leq J(H) + 2\chi(E) + d$ from Lemma 5.1. Again, we assume that all components of $U - (F \cup S)$ are $I$-bundles over an $I$-base for one of the components of $U$, since, just as we have seen, exceptional components only make the inequalities stronger. Arguing as we did in the calculation of Theorem 4.7 we have

$$2(I_v^+(V)) \leq \sum_{i=1}^{m}(J(A_i) + 2\chi(E_i) + d_i) + \sum_{i'=1}^{m}(J(B_{i'}) + 2\chi(E_{i'}) + b_{i'}) =$$

$$\sum_{i=1}^{m}(J(A_i) + (2\chi(E_i) + d_i - h_i) + h_i) + \sum_{i'=1}^{m}(J(B_{i'}) + (2\chi(E_{i'}) + b_{i'} - h_{i'}) + h_{i'}) =$$

$$\sum_{i=1}^{m}J(A_i) + \sum_{i'=1}^{m}J(B_{i'}) - 2I_\rho(U) + 2h = 4g - 4 - 2I_\rho(U) + 2h$$

Divide by 2 and add $I_\rho(U)$ to both sides and the inequality becomes

$$I_\rho(U) + I_v^+(V) \leq 2g - 2 + h$$

as required. □

## 6. Advanced computation

A further improvement in the last inequality is possible. It is based on a more sophisticated version of the argument in Lemma 5.1. A few preliminary observations are needed.

**Lemma 6.1.** *With the notation of Lemma 5.1, some annulus $A_0$ in $\xi|\partial$ is either not adjacent to a toral component or, if it is, the toral component $Y$ that it is adjacent to is longitudinal and $Y \cap (\xi|\partial) = A_0$.*

**Proof:** At least one component, $A_0 \subset \xi|\partial$, has the property that it $\partial$-compresses in the complement of $\mathcal{A}$. A $\partial$-compressing disk intersects $A_0$ in a single spanning arc, so if $A_0$ is adjacent to a toral component, that component is longitudinal. $\square$

**Lemma 6.2.** *Suppose $\Gamma$ is a bipartite graph with red and blue vertices. Then there is a subgraph $\Gamma_- \subset \Gamma$ so that*

1. *$\Gamma_-$ contains all red vertices;*
2. *each component of $\Gamma_-$ contains at most one more red vertex than blue vertices;*
3. *each blue vertex not in $\Gamma_-$ is incident to at least three components of $\Gamma_-$.*

**Proof:** The subset consisting of all red vertices satisfies the first two conditions. Let $\Gamma_-$ be a maximal subgraph satisfying these conditions. If any blue vertex is incident to one or two components of $\Gamma_-$, it could have been added without violating either condition. The result follows. $\square$

**Lemma 6.3.** *Under the hypotheses of Lemma 5.1, the last inequality can be improved to $c(Y) + s \leq J(H) + 2\chi(E) + d$.*

**Proof:** : We would like to apply the argument of Lemma 5.1 one component of $E$ at a time. A difficulty is that removing a single component $\xi^0$ may result in an annulus which is boundary parallel in the new compression body and also lies in $\xi|\partial - \xi|\partial^0$. Deleting the annulus would destroy the inductive assumption that $(\chi - \chi^0)|\partial \subset \mathcal{A}$, but leaving it in violates the hypothesis that no annulus is boundary parallel. In either case, the inductive hypothesis isn't attained.

To circumvent this difficulty, construct an abstract bipartite graph $\Gamma$ as follows: Choose a red vertex for each component of $\xi$ and a blue vertex for every longitudinal toral component which intersects $\mathcal{A}$ only in annuli that lie in $\xi|\partial$. Add an edge between a red and blue vertex for each annulus through which the corresponding components of $H - \mathcal{A}$ are adjacent. Let $\Gamma_-$ be the subgraph identified by Lemma 6.2, and let $\Gamma_0$ be the component of $\Gamma_-$ which contains the red vertex corresponding to the component of $\xi$ that contains the annulus $A_0$ from Lemma 6.1. So, on the other side of $A_0$ is either a non-toral component or a longitudinal toral component $Y$ adjacent to no other component of $E$.

Let $l$ be the number of longitudinal toral components whose boundaries intersect $\mathcal{A}$ entirely in $\xi|\partial$ (i. e. the number of blue vertices in $\Gamma$). Remove from $H$ the submanifold $W$ consisting of all components represented by vertices in $\Gamma_0$. Call the new handlebody $H'$, and examine how the numbers have changed. The change lowers $s$ by at most one more than it lowers $l$. On the other hand, it lowers $c(Y)$ no more than it lowers $d - l$. $J(H) + 2\chi(E)$ is unchanged. So, in either case, $c(Y) + s$ is lowered at most one more than

$$J(H) + 2\chi(E) + (d - l) + l = J(H) + 2\chi(E) + d.$$

In other words, (using primes to denote the relevant numbers in $H'$) we can proceed from the inductive assumption that in $H'$

$$c(Y') + s' \leq J(H') + 2\chi(E') + d'$$

only to the conclusion that in $H$

$$c(Y) + s \leq (J(H) + 2\chi(E) + d) + 1$$

and this is not quite sufficient for the inductive step.

But the following simple alteration of $H'$ can raise $c(Y')$ by one and have no effect on any of the other numbers or on the fact that $H'$ is a handlebody. Note that after removing $W$ to obtain $H'$, the annulus $A_0$ is parallel to a longitude of $H'$. This means that if a twisted solid torus is attached to $H'$ along $A_0$ the result is still a handlebody. But the operation raises $c(Y')$ by one, either because the new torus becomes a new component of $Y'$ of complexity one, or, if the component of $H' - \mathcal{A}$ to which the torus is attached is itself toroidal (hence longitudinal), it becomes twisted and so its complexity goes up by one. After this maneuver we are able to inductively conclude from

$$c(Y') + s' \leq J(H') + 2\chi(E') + d'$$

that in $H$

$$c(Y) + s \leq J(H) + 2\chi(E) + d$$

as required.

To ensure that the inductive hypothesis is intact, first disregard any (necessarily twisted) toral component which had been adjacent only to components of $W$. Any toral component that was adjacent to $W$ but still remains either is twisted, or is adjacent to a component not in $\xi$ or is adjacent to at least two other components of $\xi$ that remain. (The last corresponds to a blue vertex not in $\Gamma_-$ above.) In any case, no $\partial$-parallel annulus is adjacent to a remaining component of $\xi$, so each can be removed without changing the inductive hypothesis. □

**Theorem 6.4.** *The last inequality in 5.3 can be improved to $I_\rho(U) + I_v^\epsilon(V) \leq 2g - 2$.*

**Proof:** We argue in a manner similar to that used in 5.3, using the refinement of 6.3: $\sum_{k=1}^{t} c(Y_k) + s \leq J(H) + 2\chi(E) + d$. Again, we will assume that no component of $\partial V$ is strongly aligned and aim to show that $I_\rho(U) + I_v^+(V) \leq 2g - 2$.

For each $i, 1 \leq i \leq m$, consider $U \cap A_i$. The components can be divided into two sorts: those that are spanning (i. e. are homeomorphic to $surface \times (I : 0, 1) \subset (A_i; \partial_- A_i, \partial_+ A_i)$) and those that are non-spanning and so are homeomorphic to $(\xi, \dot\xi) \subset (A_i, \partial_+ A_i)$. Suppose here that $\xi$ lies in a non-aligned Seifert piece $U^0$ with $I$-base $E^0$. Then the $I$-base for $\xi$ is either $E^0$ or, if $\xi$ is a collar of the boundary of the $I$-bundle over $E^0$ (which is connected if $E^0$ is non-orientable), then the $I$-base for $\xi$ is the orientable double cover of $E^0$. In any case, denote by $\xi_i$ the union of all components of $U \cap A_i$ that appear in the form $(\xi, \dot\xi) \subset (A_i, \partial_+ A_i)$. Denote the $I$-base of $\xi_i$ by $E_i$. It is a union of (possibly multiple copies) of (possibly orientable double covers of) components of $E$, the $I$-base of $U$. Let $s_i = |E_i|$ and $d_i = |\partial E_i|$. For each $B_i, 1 \leq i \leq m$, similarly define $E'_i, s'_i$ and $d'_i$.

For each $i, 1 \leq i \leq m$ let $a_i$ denote the complexity of the union of all toral components of intersection of Seifert pieces with $A_i$ and similarly define $b_i$ using $B_i$.

Then Lemma 6.3 gives

$$a_i + s_i \leq J(A_i) + 2\chi(E_i) + d_i \tag{2}$$

$$b_i + s'_i \leq J(B_i) + 2\chi(E'_i) + d'_i. \tag{3}$$

Notice that if we decide to ignore one component $(\xi^0, \dot{\xi}^0) \subset (\xi_i, \dot{\xi}_i) \subset (A_i, \partial_+ A_i)$ the effect is to reduce $s_i$ by one and to lower $2\chi(E^0) + |\partial E^0|$ by at most one (since $2\chi(E^0) + |\partial E^0| \leq 1$, with equality only when $E^0$ is a pair of pants). So the two inequalities remain true even if we ignore components of $\xi_i$. Now since $M$ is closed, each component of $U$ intersects some $A_i$ in the form $(\xi^0, \dot{\xi}^0) \subset (A_i, \partial_+ A_i)$, where $\xi^0$ has the same $I$-base as $U$. Each component of $U$ intersects some $B_i$ similarly. Ignore all other components when defining $E_i$. With this alteration, we have seen that the inequality remains true, and now there is exactly one $I$-base of each component of $U$ in $\cup_{i=1}^m E_i$. Hence

$$\sum_{i=1}^m (2\chi(E_i) + 2\chi(E'_i) + d_i + d'_i - s_i - s'_i) = -2I_\rho(U).$$

Now sum both inequalities of (2) over $1 \leq i \leq m$. We have calculated in the proof of 4.7 (see equation 1) that

$$\sum_{i=1}^m c_i + \sum_{i=1}^m c'_i = 2I_v(V) + \Lambda,$$

where $V$ is the union of aligned Seifert pieces and $\Lambda$ is the total number of longitudinal torus components occuring in all $A_i \cap (V)$ or $B_i \cap (V)$ that are disjoint from $F$. We also know that

$$\sum_{i=1}^m J(A_i) + \sum_{i=1}^m J(B_i) = 4g - 4.$$

So we have

$$2I_v(V) + \Lambda + \sum_{i=1}^m (s_i + s'_i) \leq 4g - 4 + \sum_{i=1}^m (2\chi(E_i) + 2\chi(E'_i) + d_i + d'_i).$$

Now subtract from both sides $\sum_{i=1}^m (s_i + s'_i)$. Rearranging a bit we have

$$2I_v(V) + \Lambda \leq 4g - 4 + \sum_{i=1}^m (2\chi(E_i) + 2\chi(E'_i) + d_i + d'_i - s_i - s'_i) = 4g - 4 - 2I_\rho(U).$$

So we have

$$2(I_\rho(U) + I_v(V)) + \Lambda \leq 4g - 4.$$

Each component of $V$ with no exceptional fibers contributes at least two to $\Lambda$, and those with one exceptional fiber contribute at least one each to $\Lambda$. So for each component $V^0$ in $V$ that has no exceptional fibers, we can shift two from $\Lambda$ into $I_v(V^0)$ to get $I_v^+(V^0)$, and for each component with one exceptional fiber, we can shift one. Drop all the rest of $\Lambda$ from the inequality to get

$$2(I_\rho(U) + I_v^+(V)) \leq 4g - 4$$

as required. □

**Definition 6.5.** *A toroidal Seifert piece is called* weakly toroidal *if each incompressible torus is either boundary parallel or bounds an I-bundle over the Klein bottle.*

A weakly toroidal piece either
- has base the once punctured Möbius band and has no exceptional fibers or
- has base the Möbius band and has one exceptional fiber.

In particular, a Seifert piece $V^0$ is weakly toroidal if and only if $I_v^\epsilon(V^0) = 3/2$

**Corollary 6.6.** *Suppose $M$ is a closed orientable irreducible 3-manifold that has a genus $g$ Heegaard splitting. Let $\Theta$ be the decomposing collection of tori for $M$. Suppose there are $n'$ weakly toroidal Seifert pieces of $M - \Theta$ and $n$ other toroidal pieces.*

*Then the number of non-Seifert components of $M - \Theta$ is at most $g - 1$ and the total number of components is at most $3g - 3 - n - n'/2$.*

**Proof:** : The first inequality of 5.3 implies that $p + h \leq g - 1$. Now examine the inequality of 6.4. Note first that if $V^0$ is a Seifert component with $I_v^\epsilon(V^0) = 1$, then $V^0$ is atoroidal. So if $v$ of the $r$ components of $V$ are toroidal but not weakly toroidal and $n'$ are weakly toroidal, then the inequality gives

$$I_\rho(U) + r \leq 2g - 2 - v - n'/2.$$

Now for each component $U_0$ of $U$, $I_\rho(U_0) \geq 1$ unless $U_0$ has $I$-base a pair of pants (and so is atoroidal). So if $h'$ of the components of $U$ are toroidal we have

$$r \leq 2g - 2 - v - n'/2 - h'.$$

Combining, we get

$$p + h + r \leq 3g - 3 - v - n'/2 - h' = 3g - 3 - n - n'/2,$$

as required. □

A defect of Theorems 5.3 and 6.4 is that one doesn't know *ab initio* whether a Seifert piece of $M - \Theta$ will be aligned or not. So it would be useful to have as broad a statement as possible which does not require knowing which pieces are aligned. Corollary 6.6 is an example, but it can be generalized.

**Lemma 6.7.** *Suppose $W$ is a connected orientable Seifert manifold with non-empty boundary and suppose $W$ is an I-bundle construct over $E$, with $\chi(E) \leq -1$. Let $\rho$ denote the slope of $\partial E$ in $\partial W$. Then either $W$ is a product bundle and $E$ is a pair of pants or*

$$I_v^+(W) \leq I_\rho(W) - \chi(E).$$

**Proof:** Recall the definitions: $I_\rho(W) = 1 - 2\chi(E) - d$ where $d$ is the number of components in $\partial E$. If the base space of $W$ is $P$ and there are $f$ exceptional fibers, then $I_v^+(W)$ is either $f - \chi(P)$ if there is more than one exceptional fibers, $1/2 - \chi(P)$ if there is exactly one and $1 - \chi(P)$ if there are no exceptional fibers. So, as an alternative formulation of the inequality, it suffices to show that

$$1 - d - f + \chi(P) - 3\chi(E) \geq 0$$

if $f \geq 1$ and
$$-d + \chi(P) - 3\chi(E) \geq 0$$
when $f = 0$.

First suppose there are no exceptional fibers, so $E$ is a cover of $P$ of degree $\mu$, say. Furthermore $\chi(E) \leq -1$, so $\chi(P) \leq -1$, and $-d - \chi(E) = -\chi(\hat{E})$, where $\hat{E}$ is the closed surface obtained by capping off all boundary components of $E$. So in order to show that
$$-d + \chi(P) - 3\chi(E) \geq 0,$$
it suffices to show that

(4) $$\chi(P) - 2\mu\chi(P) = (1 - 2\mu)\chi(P) \geq \chi(\hat{E}).$$

Noting that $\chi(\hat{E}) \leq 2$, this inequality follows immediately unless $\mu = 1$ (so $P = E$), and even then it is immediate unless $\chi(P) = -1$. When $\mu = 1$ and $\chi(P) = -1$ there are three possibilities: If $P$ is a once-punctured torus or Klein bottle, then $\chi(\hat{E}) = 0$ and the inequality (4) still follows. Similarly, if $P$ is a twice-punctured projective plane, then $\chi(\hat{E}) = 1$ and the inequality (4) again follows. This leaves only the final case, when $P = E$ is a pair of pants and $W$ is a product bundle, which is a case allowed by Lemma 6.7.

Now consider the case in which there are exceptional fibers. $E$ is a manifold cover of an orbifold whose underlying surface is the base $P$. The Riemann-Hurwitz formula describes the relation between the Euler characteristics ([Sco, p. 427]) of $P$ and $E$. Suppose that the covering map $E \to P$ is of degree $\mu$, and that the exceptional fibers are of type $p_i/q_i, 1 \leq i \leq f$. Then

$$\chi(E) = \mu(\chi(P) - \sum_{i=1}^{f}(1 - 1/q_i)).$$

So we are asked to verify that

$$1 - d - f + \chi(P) - 3\mu(\chi(P) - \sum_{i=1}^{f}(1 - 1/q_i)) \geq 0.$$

First notice that this expression is minimized when each $q_i$ is as small as possible so in order to prove the inequality it suffices to set $q_i = 2$, i. e. $\sum_{i=1}^{f}(1 - 1/q_i) = f/2$. So it suffices to prove

$$1 - d - f + \chi(P) - 3\mu(\chi(P) - f/2) \geq 0.$$

We distinguish three cases:

**Case 1:** $\chi(P) = 0$. Then $P$ is an annulus or a Möbius band. If $f = 1$ then we can calculate directly:
$$I_\rho(W) - \chi(E) \geq -\chi(E) \geq 1 = I_v^+(W)$$
as required. So we may as well assume that $f \geq 2$.

$P$ has one or two boundary components. This implies that $d \leq \mu|\partial P| \leq 2\mu$, so, using also $\mu \geq 2$ we have

$$\begin{aligned}
1 - d - f + \chi(P) - 3\mu(\chi(P) - f/2) &\geq 1 - 2\mu - f + 3\mu f/2 = \\
1 - 2\mu + (3\mu/2 - 1)f &\geq 1 - 2\mu + (3\mu/2 - \mu/2)f = \\
1 - 2\mu + \mu f &= 1 + (f - 2)\mu \geq 0
\end{aligned}$$

as required.

**Case 2:** $\chi(P) \leq -1$. We don't know how many components $\partial P$ has, but we do know that $d \leq \mu|\partial P|$. Then it suffices to prove that

$$1 - \mu|\partial P| - f + \chi(P) - 3\mu\chi(P) + 3\mu f/2 \geq 0.$$

The left expression can be rearranged to become

$$\begin{aligned}
1 - \mu(|\partial P| + \chi(P)) + \chi(P) - 2\mu\chi(P) + (3\mu/2 - 1)f &\geq \\
1 - \mu\chi(\hat{P}) + (1 - 2\mu)\chi(P) + (3\mu/2 - \mu/2)f &\geq \\
(1 - 2\mu)(1 + \chi(P)) + \mu f &\geq 0
\end{aligned}$$

as required.

**Case 3:** $\chi(P) = 1$. Then $P$ is a disk, so $f \geq 2$. If $f = 2$ then we can calculate directly:

$$I_\rho(W) - \chi(E) \geq -\chi(E) \geq 1 = I_v^+(W)$$

as required. So we may as well assume that $f \geq 3$.

Since $P$ is a disk, it has one boundary component. This implies that $d \leq \mu$, so

$$1 - d - f + \chi(P) - 3\mu(\chi(P) - f/2) \geq 2 - \mu - f - 3\mu(1 - f/2) = 2 - f - \mu(4 - 3f/2).$$

Since $f \geq 3, 4 - 3f/2 \leq 0$ and the expression $2 - f - \mu(4 - 3f/2)$ is minimized when $\mu = 2$. In that case we have

$$2 - f - \mu(4 - 3f/2) \geq -6 + 2f \geq 0,$$

as required. □

**Corollary 6.8.** *Suppose $M$ is a closed orientable irreducible 3-manifold that has a genus $g$ Heegaard splitting. Let $\Theta$ be the decomposing collection of tori for $M$. Suppose $p$ of the components of $M - \Theta$ are non-Seifert, and $s$ are of the form $S_3 \times S^1$, where $S_3$ is the three-punctured sphere. Let $W$ denote the union of the remaining Seifert pieces. Then*

$$p + s + I_v^\epsilon(W) \leq 3g - 3.$$

**Proof:** As in the proofs of Theorems 4.7, 5.3, and 6.4, untelescope a minimal genus Heegaard splitting for $M$. Let $V$ be the union of aligned Seifert pieces of $M - \Theta$. Let $s'$ be the number of non-aligned components that are of the form $S_3 \times S^1$, let $U$ denote the union of the remaining non-aligned Seifert components, and let $E$ be the $I$-base of $U$. (Recall from the note before 3.7 that we can assume each component of $E$ has negative Euler characteristic.) Then from Theorem 5.3 we have

$$p + s' - \chi(E) \leq g - 1.$$

From Theorem 6.4 we have
$$I_\rho(U) + I_v^\epsilon(V) \leq 2g - 2.$$

Apply Lemma 6.7 and add the inequalities to get
$$p + s' + I_v^+(U) + I_v^\epsilon(V) \leq p + s' - \chi(E) + I_\rho(U) + I_v^\epsilon(V) \leq 3g - 3.$$

Now let $W$ be the union of $U$ and all components of $V$ which are not of the form $S_3 \times S^1$, and suppose $s''$ of the components of $V$ are of the form $S_3 \times S^1$. Then, since $I_v^\epsilon(U) \leq I_v^+(V)$ the previous inequality can be written

$$p + s' + 2s'' + I_v^\epsilon(W) \leq 3g - 3,$$

and this certainly implies that

$$p + s + I_v^\epsilon(W) = p + s' + s'' + I_v^\epsilon(W) \leq 3g - 3$$

as required.  □

Note that this generalizes Corollary 6.6, since any toroidal Seifert piece has epsilon vertical index $\geq 2$, except the weakly toroidal pieces with epsilon vertical index $= 3/2$.

**Corollary 6.9.** *Suppose $M$ is a closed orientable irreducible 3-manifold that has a genus $g$ Heegaard splitting. Let $\Theta$ be the decomposing collection of tori for $M$. Suppose $p$ of the components of $M - \Theta$ are not Seifert fibered spaces, and the union of Seifert components of $M - \Theta$ has base space $P$ and $f$ exceptional fibers. Then*
$$f - \chi(P) \leq 3g - 3 - p.$$

**Proof:** Apply Corollary 6.8: $-\chi(S_3) = 1$, so
$$f - \chi(P) = s + I_v(W) \leq s + I_v^\epsilon(W) \leq 3g - 3 - p.$$

□

MARTIN SCHARLEMANN, MATHEMATICS DEPARTMENT, UNIVERSITY OF CALIFORNIA, SANTA BARBARA, CA USA
  *E-mail address*: `mgscharl@math.ucsb.edu`

JENNIFER SCHULTENS, MATHEMATICS DEPARTMENT, EMORY UNIVERSITY, ATLANTA, GA 30322 USA
  *E-mail address*: `jcs@mathcs.emory.edu`